\documentstyle[12pt]{article}
\newcommand{\sect}[1]{\section{#1}\setcounter{equation}{0}}

\font\mbn=msbm10 scaled \magstep1
\font\mbs=msbm7 scaled \magstep1
\font\mbss=msbm5 scaled \magstep1
\newfam\mbff
\textfont\mbff=\mbn
\scriptfont\mbff=\mbs
\scriptscriptfont\mbff=\mbss\def\mbf{\fam\mbff}

\def\Co{{\mbf C}}
\def\To{{\mbf T}}

\def\N{{\mbf N}}

\newtheorem{Th}{Theorem}[section]
\newtheorem{Lm}[Th]{Lemma}

\newtheorem{Proposition}[Th]{Proposition}
\newtheorem{R}[Th]{Remark}

\author{Alexander Brudnyi\thanks{Research supported in part 
by Max-Planck-Institut f\"{u}r Mathematik.
}\\
Department of Mathematics and Statistics\\
University of Calgary, Calgary\\
Canada}
\title{Hartogs Type Theorems for $CR$ $L^{2}$ functions
on Coverings of Strongly Pseudoconvex Manifolds}
\date{} 
\begin{document} 
\maketitle
\begin{abstract}
{We prove an analog of the classical Hartogs extension theorem for $CR$
$L^{2}$ functions defined on boundaries of
certain (possibly unbounded) domains on coverings of strongly
pseudoconvex manifolds. Our result is related to a question formulated
in the paper of Gromov, Henkin and Shubin [GHS] on holomorphic
$L^{2}$ functions on coverings of pseudoconvex manifolds.
}
\end{abstract}

\sect{\hspace*{-1em}. Introduction.}
{\bf 1.1.} In this paper, 
following our previous work [Br4], 
we continue to study holomorphic $L^{2}$ functions on coverings of strongly pseudoconvex manifolds. The subject was originally motivated by the paper [GHS] of Gromov, Henkin and Shubin.
In [GHS] the von Neumann dimension was used to measure the space of holomorphic $L^{2}$ functions on {\em regular} (i.e., {\em Galois}) coverings of a strongly pseudoconvex manifold $M$. In particular, it was shown that the space of such functions is infinite-dimensional. It was also asked whether the regularity of the covering is relevant for the existence of many holomorphic $L^{2}$ functions on $M'$ or it is just an artifact of the chosen method of
the proof which requires a use of von Neumann algebras. 

In an earlier paper [Br4] we proved that actually the regularity of $M'$ is irrelevant for the existence of many holomorphic $L^{2}$ functions on $M'$. Moreover, we obtained an extension of some of the main results of [GHS]. The method of the proof used in [Br4] is completely different and (probably) easier than that used in [GHS] and is based on $L^{2}$ cohomology techniques, as well as, on the geometric properties of $M$. Also, in [Br1]-[Br3] the case of coverings of pseudoconvex domains in Stein manifolds was considered.  Using the methods of the theory of coherent Banach sheaves together with Cartan's vanishing cohomology theorems, we proved some results on holomorphic $L^{p}$ functions, $1\leq p\leq\infty$, defined on such coverings.\\
{\bf 1.2.} The present paper is related to one of the open problems posed in [GHS], a Hartogs type theorem for coverings of strongly pseudoconvex manifolds. Let us recall that for a bounded open set
$D\subset\Co^{n}\ (n>1)$ with a connected smooth boundary
$bD$ the classical Hartogs theorem states that any holomorphic function
in some neighbourhood of $bD$ can be extended to a holomorphic function on
a neighbourhood of the closure $\overline{D}$. 
In [Bo] Bochner proved a similar extension result for functions defined 
on the $bD$ only. In modern language his result says that for a smooth 
function
defined on the $bD$ and satisfying the tangential Cauchy-Riemann equations
there is an extension to a holomorphic function in $D$ which is smooth
on $\overline{D}$. In fact, this statement follows from 
Bochner's proof (under some smoothness conditions). However at that time
there was not yet the notion of a $CR$-function. Over the years significant
contributions to the area of Hartogs theorem were made by many prominent
mathematicians, see the history and the references in the paper of
Harvey and Lawson [HL, Section 5]. A general Hartogs-Bochner type
theorem for bounded domains $D$ in Stein manifolds was proved by Harvey and
Lawson [HL, Theorem 5.1]. The proof relies heavily upon the fact that
for $n\geq 2$ any $\overline{\partial}$-equation with compact support on
an $n$-dimensional Stein manifold has a compactly supported solution.
In [Br2] and [Br3] we proved some extensions of the theorem of Harvey and
Lawson for certain (possibly unbounded) domains on coverings of Stein 
manifolds. In the present paper we prove an analogous result for $CR$ 
$L^{2}$ functions defined on boundaries
of certain domains on coverings of strongly 
pseudoconvex manifolds. More general Hartogs type theorems
for $CR$-functions of slow growth on boundaries of such domains will be
presented in a forthcoming paper.\\
{\bf 1.3.} 
Let $M\subset\subset N$ be a domain with smooth boundary 
$bM$ in an $n$-dimensional complex manifold $N$, specifically,
\begin{equation}\label{m1}
M=\{z\in N\ :\ \rho(z)<0\}
\end{equation}
where $\rho$ is a real-valued function of class $C^{2}(\Omega)$ in a
neighbourhood $\Omega$ of the compact set $\overline{M}:=M\cup bM$
such that
\begin{equation}\label{m2}
d\rho(z)\neq 0\ \ \ {\rm for\ all}\ \ \ z\in bM\ .
\end{equation}
Let $z_{1},\dots, z_{n}$ be complex local coordinates in $N$ near $z\in bM$.
Then the tangent space $T_{z}N$ at $z$ is identified with $\Co^{n}$.
By $T_{z}^{c}(bM)\subset T_{z}N$ we denote the complex tangent space to
$bM$ at $z$, i.e.,
\begin{equation}\label{m3}
T_{z}^{c}(bM)=\{w=(w_{1},\dots,w_{n})\in T_{z}(N)\ :\ \sum_{j=1}^{n}
\frac{\partial\rho}{\partial z_{j}}(z)w_{j}=0\}\ .
\end{equation}
The {\em Levi form} of $\rho$ at $z\in bM$ is a hermitian form on
$T_{z}^{c}(bM)$ defined in local coordinates by the formula
\begin{equation}\label{m4}
L_{z}(w,\overline{w})=\sum_{j,k=1}^{n}
\frac{\partial^{2}\rho}{\partial z_{j}\partial\overline{z}_{k}}(z)w_{j}
\overline{w}_{k}\ .
\end{equation}
The manifold $M$ is called {\em pseudoconvex} if $L_{z}(w,\overline{w})\geq 0$
for all $z\in bM$ and $w\in T_{z}^{c}(bM)$. It is called {\em strongly
pseudoconvex} if $L_{z}(w,\overline{w})>0$ for all $z\in bM$ and all
$w\neq 0$, $w\in T_{z}^{c}(bM)$.

Equivalently, strongly pseudoconvex manifolds can be described as the ones 
which locally, in a neighbourhood of any boundary point, can be presented as 
strictly convex domains in $\Co^{n}$. It is also known (see [C], [R]) that 
any strongly pseudoconvex manifold admits a proper holomorphic map with 
connected fibres onto a normal Stein space. In particular, if
$M$ is a strongly pseudoconvex non-Stein manifold of
complex dimension $n\geq 2$, then the union $C_{M}$
of all compact complex subvarieties of $M$
of complex dimension $\geq 1$ is a compact complex subvariety
of $M$. 

Let $r: M'\to M$ be an unbranched covering of $M$.
Assume that $N$ is equipped with a
Riemannian metric $g_{N}$. By $d$ we denote the path metric on $M'$ 
induced by the pullback of $g_{N}$. Consider a domain 
$\widetilde D\subset\subset M$ 
with a connected $C^{1}$ smooth boundary $b\widetilde D$ such that
\begin{equation}\label{inter}
b\widetilde D\cap C_{M}=\emptyset.
\end{equation}
Let $D$ be a connected component of $r^{-1}(\widetilde D)$. 
By $bD$ we denote the 
boundary of $D$ and by $\overline{D}\subset M'$ the closure of $D$.
Also, by ${\cal O}(D)$ we denote the space of holomorphic functions
on $D$. Now, recall that a continuous
function $f$ on $bD$ is called $CR$ if for every smooth $(n,n-2)$-form
$\omega$ on $M'$ with compact support one has
$$
\int_{bD}f\cdot\overline{\partial}\omega=0\ \! .
$$
If $f$ is smooth this is equivalent to $f$
being a solution of the tangential $CR$-equations: 
$\overline{\partial}_{b}f=0$ (see, e.g., [KR]).

Let $dV_{M'}$ and $dV_{bD}$ be the Riemannian volume forms on $M'$ and $bD$ obtained by the 
pullback of the Riemannian metric $g_{N}$ on $N$.
By $H^{2}(D)$ we denote the Hilbert space of holomorphic functions $g$ 
on $D$ with norm
$$
\left(\int_{z\in D}|g(z)|^{2}dV_{M'}(z)\right)^{1/2}\ .
$$
Also, $L^{2}(bD)$ stands for the Hilbert space of functions $g$ on 
$bD$ with norm
$$
\left(\int_{z\in bD}|g(z)|^{2}dV_{bD}(z)\right)^{1/2}\ .
$$
The following question was asked in [GHS, Section 4]:

{\em Suppose that $D$ is a regular covering of a strongly pseudoconvex manifold
$\widetilde D\subset\subset M$. Is it true that
for every $CR$-function $f\in L^{2}(bD)\cap C(\overline{D})$ there exists
$F\in H^{2}(D)\cap C(\overline{D})$ such that $F|_{bD}=f$?}

In the present paper we give a particular answer to this question.
To formulate our results we require the following definitions. 

For every
$x$ from the closure of $\widetilde D$ we introduce the Hilbert space $l_{2,x}(D)$ of functions 
$g$ on $r^{-1}(x)\cap\overline{D}$ with norm
\begin{equation}\label{norm}
|g|_{x}:=\left(\sum_{y\in
r^{-1}(x)\cap\overline{D}}|g(y)|^{2}\right)^{1/2} .
\end{equation}
Next, we 
introduce the Banach space ${\cal H}_{2}(D)$ of 
holomorphic on $D$ functions $f$ with norm
$$
|f|_{D}:=\sup_{x\in\widetilde D}|f|_{x}\ .
$$
Similarly, we introduce the Banach space ${\cal L}_{2}(bD)$ of 
continuous on $bD$ functions $g$ with norm
$$
|g|_{bD}:=\sup_{x\in b\widetilde D}|f|_{x}.
$$

Let ${\cal U}=(U_{i})_{i\in I}$ be a finite open cover of
$b\widetilde D$ by open simply connected sets $U_{i}\subset\subset M$. 
Then $r^{-1}(U_{i})\cap bD$ is homeomorphic to 
$(U_{i}\cap b\widetilde D)\times Q$ where $Q$ is the fibre of the covering
$r:D\to\widetilde D$. In what follows we identify $r^{-1}(U_{i})\cap bD$
with $(U_{i}\cap b\widetilde D)\times Q$.

Suppose that $f\in C(bD)$ is a $CR$-function satisfying the
following conditions
\begin{itemize}
\item[(1)]
$$
f\in {\cal L}_{2}(bD);
$$
\item[(2)]
for any $i\in I$ and any $z_{1},z_{2}\in b\widetilde D\cap U_{i}$ 
there is a constant $L_{i}$ such that
$$
\left(\sum_{q\in Q}\left|\frac{f(z_{1},q)-f(z_{2},q)}{d((z_{1},q),(z_{2},q))}
\right|^{2}\right)^{1/2}\leq 
L_{i}.
$$
\end{itemize}
(It is easy to show that condition (2) is independent of the choice of the
cover.)
\begin{Th}\label{te1}
For any $CR$-function $f$ on $bD$ satisfying conditions (1) and (2) 
there exists $\hat f\in {\cal H}_{2}(D)\cap C(\overline{D})$ 
such that
$$
\hat f|_{bD}=f\ \ \ {\rm and}\ \ \ |\hat f|_{D}=|f|_{bD}. 
$$
\end{Th}
\begin{R}\label{re1}
{\em (A) If, in addition, $bD$ is smooth of class $C^{k}$, 
$1\leq k\leq\infty$, and $f\in C^{s}(bD)$, $1\leq s\leq k$, then the 
extension $\hat f$
belongs to ${\cal O}(D)\cap C^{s}(\overline{D})$. This follows from
[HL, Theorem 5.1].\\
(B) From the Cauchy
integral formula it follows that the hypotheses of the theorem
are true if $f$
is the restriction to $bD$ of a holomorphic function from 
${\cal H}_{2}(W)$ where $\widetilde W:=r(W)\subset\subset M$ is a 
neighbourhood of $b\widetilde D$ and $W$ is a connected component of
$r^{-1}(\widetilde W)$ containing $bD$ (see [Br1, Proposition 2.4] for
similar arguments).\\
(C) It was shown in [Br4, Theorem 1.1] that holomorphic functions 
from ${\cal H}_{2}(M')$
separate points on $M'\setminus C_{M}'$ where
$C_{M}':=r^{-1}(C_{M})$. Thus there are sufficiently
many $CR$-functions $f$ on $bD$ satisfying conditions (1) and (2).}
\end{R}

As before by ${\cal L}_{2}(M')$ we denote the Banach space of continuous
functions $f$ on $M'$ with norm
$$
|f|_{M'}:=\sup_{x\in M}|f|_{x}.
$$
where $|\cdot|_{x}$, $x\in M$, is defined as in (\ref{norm}) with $M'$ 
substituted for $\overline{D}$.
Also, for a measurable locally bounded 
$(0,1)$-differential form $\eta$ on $M'$ by $|\eta|_{z}$,
$z\in M'$, we denote the norm of $\eta$ at $z$ defined by the natural
hermitian metric on the fibres of the cotangent bundle $T^{*}M'$ on
$M'$. We say that such $\eta$ belongs to the space ${\cal E}_{2}(M')$ if 
\begin{equation}\label{cond2}
|\eta|_{M'}:=\sup_{x\in M}\left(\sum_{z\in r^{-1}(x)}
|\eta|_{z}^{2}\right)^{1/2}<\infty .
\end{equation}
(Note that this definition does not depend on the choice of the Riemannian metric on $N$, and that the expression in the brackets is 
correctly defined for almost all $x\in M$.)\\ 
By $supp\ \eta$ we denote support of 
$\eta$, i.e., the minimal closed set $K\subset M'$ such that $\eta$ equals
zero almost everywhere on $M'\setminus K$.

As mentioned above, the proof of the classical Hartogs theorem is based on
the fact that for $n\geq 2$ any $\overline\partial$-equation with compact
support on an $n$-dimensional Stein manifold has a compactly supported
solution. Similarly our proof of Theorem \ref{te1} is based on the 
following result.
\begin{Th}\label{te2}
Let $O\subset\subset M\setminus C_{M}$. 
Assume that a $(0,1)$-form $\eta$ on $M'$ belongs to
${\cal E}_{2}(M')$, is $\overline\partial$-closed 
(in the distributional sense) and  
$$
r(supp\ \eta)\subset O.
$$
Then there are a function $F\in {\cal L}_{2}(M')$ and a neighborhood $U\subset M$ of $bM$
such that $\overline\partial F=\eta$ (in the distributional sense)
and $F|_{r^{-1}(U)}=0$.
\end{Th}
(Since $M'$ can be thought of as a subset of a covering $L'$ of a neighbourhood $L$ of $\overline{M}$, the boundary $bM'$
of $M'$ is correctly defined.)
\begin{R}\label{re2}
{\rm (A) Condition (2) in the formulation
of Theorem \ref{te1} means that $f$ is a Lipschitz section of a Hilbert vector
bundle on $b\widetilde D$ with fibre $l_{2}(Q)$ associated with the natural
action of the fundamental group $\pi_{1}(b\widetilde D)$ of $b\widetilde D$ on
$l_{2}(Q)$ (see [Br1, Example 2.2(b)] for a similar construction). This
condition is required by the method of the proof. It would be
interesting to know to what extent it is necessary.\\
(B) Another interesting question is whether a general extension theorem for 
$CR$-functions on $bD$ without growth condition might hold.}
\end{R}

{\em Acknowledgment.} This work was written during my stay at the Max-Planck-Institut f\"{u}r Mathematik in Bonn. I am deeply grateful to MPIM for
hospitality and financial support.
\sect{\hspace*{-1em}. Proof of Theorem \ref{te1}.}
In this section we prove Theorem \ref{te1} modulo Theorem \ref{te2}. Then
in the next section we prove Theorem \ref{te2}.

Since $b\widetilde D$ is a compact $C^{1}$ smooth manifold, there are a 
neighbourhood $O\subset\subset M\setminus C_{M}$ of 
$b\widetilde D$ and a $C^{1}$ retraction
$p:O\to b\widetilde D$. (As such $O$ one can take, e.g., a neighbourhood
of the zero section of the normal vector bundle on $b\widetilde D$.) Without
loss of generality we may assume also that fundamental groups 
$\pi_{1}(O)$ and $\pi_{1}(b\widetilde D)$ are isomorphic. Let
$O'$ be a connected component of $r^{-1}(O)\subset M'$ containing $bD$.
Then by the covering homotopy theorem there is a $C^{1}$ retraction
$p':O'\to bD$ such that $r\circ p'=p\circ r$.

Let $\rho$, $0\leq\rho\leq 1$, 
be a $C^{\infty}$ function on $M$ equals 1 in a neighbourhood of 
$b\widetilde D$ with $supp\ \!\rho\subset\subset O$. Consider the
$C^{\infty}$ function $\rho':=\rho\circ r$ on $M'$.

Let ${\cal V}=(V_{j})_{j\in J}$ be a finite open cover of 
$\widetilde D\cup b\widetilde D$ by simply connected
coordinate charts $V_{j}\subset\subset M$. We naturally identify 
$r^{-1}(V_{j})$ with $V_{j}\times S$ where $S$ is the fibre of 
$r:M'\to M$. Then in these local coordinates on $M'$ we have
\begin{equation}\label{e2.1}
p'(z,s)=(p(z),s),\ \ \ \rho'(z,s)=\rho(z),\ \ \ (z,s)\in O'\cap r^{-1}(V_{j}),\ \
j\in J.
\end{equation}

Next, for a $CR$-function $f$ satisfying the assumptions of the theorem
we define
\begin{equation}\label{e2.2}
f_{1}(z):=\rho'(z)\cdot f(p'(z)),\ \ \ z\in\overline{D}.
\end{equation}
\begin{Lm}\label{le2.1}
In the above local coordinates on $M'$ one has
$$
\left(\sum_{s\in S}\left|\frac{f_{1}(z_{1},s)-f_{1}(z_{2},s)}
{d((z_{1},s),(z_{2},s))}
\right|^{2}\right)^{1/2}\leq 
C_{j},\ \ \ (z_{1},s),\ (z_{2},s)\in\overline D\cap r^{-1}(V_{j}),\ \
j\in J,
$$
for some numerical constants $C_{j}$.
\end{Lm}
{\bf Proof.}
By $d_{N}$ we denote the path metric on $N$ determined by the Riemannian metric
$g_{N}$. Since the path metric $d$ on $M'$ is obtained by the pullback of
$g_{N}$, we have $d((z_{1},s),(z_{2},s))=d_{N}(z_{1},z_{2})$.
Also, by the definition of $p'$ and $\rho'$ we clearly have for some $C>0$,
$$
\begin{array}{c}
\displaystyle
d(p'(z_{1},s),p'(z_{2},s))\leq C d_{N}(z_{1},z_{2})\ \ \ {\rm for\ all}
\ \ \ z_{1}, z_{2}\in supp\ \rho,\ \ \ {\rm and}\\
\\
\displaystyle
|\rho'(z_{1},s)-\rho'(z_{2},s)|\leq C d_{N}(z_{1},z_{2})\ \ \
{\rm for\ all}\ \ \ z_{1}, z_{2}\in M.
\end{array}
$$

Using these inequalities, condition (2) of the theorem and
the triangle inequality for $l_{2}$ norms we obtain that there is $A>0$ such
that for $z_{1}, z_{2}\in supp\ \!\rho$
$$
\begin{array}{c}
\displaystyle \left(\sum_{s\in S}
\left|\frac{f_{1}(z_{1},s)-f_{1}(z_{2},s)}{d_{N}(z_{1},z_{2})}
\right|^{2}
\right)^{1/2}\leq\\
\\
\displaystyle
\left(\sum_{s\in S}\left\{\left|\frac{\rho(z_{1})-\rho(z_{2})}{d_{N}(z_{1},z_{2})}
\right|\cdot |f(p(z_{1}),s)|+
|\rho(z_{2})|\cdot\left|\frac{f(p(z_{1}),s)-f(p(z_{2}),s)}
{d_{N}(z_{1},z_{2})}\right|\right\}^{2}\right)^{1/2}\leq\\
\\
\displaystyle C\left\{\left(\sum_{s\in S}|f(p(z_{1}),s)|^{2}
\right)^{1/2}+\left(\sum_{s\in S}\left|
\frac{f(p(z_{1}),s)-f(p(z_{2}),s)}
{d((p(z_{1}),s),(p(z_{2}),s))}
\right|^{2}\right)^{1/2}\right\}\leq A.
\end{array}
$$

Suppose now that, e.g., $z_{1}\in supp\ \rho$ and $z_{2}\not\in supp\ \rho$.
Then the term with $|\rho(z_{2})|$ in the second line of
the above inequalities disappears and again we get the require estimate.
Finally, the case $z_{1},z_{2}\not\in supp\ \rho$ is obvious.\ \ \ \ \ $\Box$

This lemma in particular implies that $f_{1}$ is a bounded Lipschitz function
on $\overline{D}$. Now, using the McShane extension theorem [M] we extend
$f_{1}$ to a Lipschitz function $\widetilde f$ on $M'$.

Further, since locally the metric $d$ is equivalent to the Euclidean metric and
since $\widetilde f$ is Lipschitz 
on $M'$, by the Rademacher theorem, see, e.g., [Fe, Section 3.1.6],
$\widetilde f$ is differentiable almost everywhere.
In particular, $\overline\partial\widetilde f$ is a 
$(0,1)$-form on $M'$ whose coefficients in its 
local coordinate representations 
are $L^{\infty}$-functions. Let $\chi_{D}$ be the characteristic function of
$D$. Consider the $(0,1)$-form on $M'$ defined by 
$$
\omega:=\chi_{D}\cdot\overline\partial\widetilde f\ .
$$
Then repeating word-for-word the arguments of [Br3, Lemma 3.3] we get
\begin{Lm}\label{le2}
$\omega$ is $\overline\partial$-closed in the distributional sense.\ \ \ \ \ 
$\Box$
\end{Lm}

Also, the inequality of Lemma \ref{le2.1} implies immediately that
$\omega\in {\cal E}_{2}(M')$, see (\ref{cond2}). Moreover, by our construction
$r(supp\ \!\omega)\subset\subset M\setminus C_{M}$.
 Thus according to
Theorem \ref{te2} there is a continuous function $F\in {\cal L}_{2}(M')$ such
that $\overline\partial F=\omega$ and $F|_{r^{-1}(U)}=0$ for a neighbourhood
$U\subset M$ of $bM$. Since $D\subset M'$ is a domain with a connected
boundary, and $F$ is holomorphic outside $\overline{D}$ (by the definition
of $\omega$), the latter
implies that $F|_{bD}=0$.

We set
$$
\hat f(z):=f_{1}(z)- F(z),\ \ \ z\in\overline{D}.
$$
Using the above properties of $f_{1}$ and $F$ one obtains easily that
$$
\hat f\in {\cal O}(D)\cap C(\overline{D})\ \ \ {\rm and}\ \ \ \hat f|_{bD}=f.
$$
Since $f_{1}$ and $F|_{\overline{D}}$ belong to ${\cal L}_{2}(\overline{D})$,
$\hat f\in {\cal H}_{2}(D)$.
Now, the identity $|\hat f|_{D}=|f|_{bD}$ follows from the fact that
the function $z\mapsto |f|_{z}$, $z\in \widetilde D\cup b\widetilde D$,
see (\ref{norm}), is
continuous and plurisubharmonic on $\widetilde D$.

This completes the proof of the theorem.\ \ \ \ \ $\Box$
\sect{\hspace*{-1em}. Proof of Theorem \ref{te2}.}
{\bf 3.1.} In Sections 3.1-3.6 we collect some auxiliary results required in 
the proof. Then in Section 3.7 we prove the theorem.
 
Let $X$ be a complete K\"{a}hler manifold of dimension $n$ with a 
K\"{a}hler form $\omega$ and
$E$ be a hermitian holomorphic vector bundle on $X$ with curvature $\Theta$. 
Let $L_{2}^{p,q}(X,E)$ be the space of $L^{2}$ $E$-valued $(p,q)$-forms on $X$
with the $L^{2}$ norm, and let $W_{2}^{p,q}(X,E)$ be the subspace of forms
such that $\overline\partial\eta$ is $L^{2}$. (The forms $\eta$ may be taken 
to be either smooth or just measurable, in which case 
$\overline\partial\eta$ is understood in the distributional sense.) The
cohomology of the resulting $L^{2}$ Dolbeault complex 
$(W_{2}^{\cdot ,\cdot},\overline\partial)$ is the $L^{2}$ cohomology
$$
H_{(2)}^{p,q}(X,E)=Z_{2}^{p,q}(X,E)/B_{2}^{p,q}(X,E)\ ,
$$
where $Z_{2}^{p,q}(X,E)$ and $B_{2}^{p,q}(X,E)$ are the spaces of 
$\overline\partial$-closed and $\overline\partial$-exact forms in
$L_{2}^{p,q}(X,E)$, respectively.

If $\Theta\geq \epsilon\omega$
for some $\epsilon>0$ in the sense of Nakano, then 
the $L^{2}$ Kodaira-Nakano vanishing theorem, see [D], [O], 
states that
\begin{equation}\label{e1}
H_{(2)}^{n,r}(X,E)=0\ \ \ {\rm for}\ \ \ r>0\ .
\end{equation}

Assume now that $\Theta\leq -\epsilon\omega$ for some $\epsilon>0$ in the
sense of Nakano. Then using a duality argument and the Kodaira-Nakano vanishing theorem (\ref{e1}) one obtains, see [L, Corollary 2.4], 
\begin{equation}\label{e1'}
H_{(2)}^{0,r}(X,E)=0\ \ \ {\rm for}\ \ \ r<n.
\end{equation} 
{\bf 3.2.}
Let $M\subset\subset N$ be a strongly pseudoconvex manifold. Without
loss of generality we will assume that $\pi_{1}(M)=\pi_{1}(N)$ and
$N$ is strongly pseudoconvex, as well. Then there exist a normal
Stein space $X_{N}$, a proper holomorphic surjective map $p:N\to X_{N}$ with
connected fibres and points $x_{1},\dots, x_{l}\in X_{N}$ such that
$$
p:N\setminus\bigcup_{1\leq i\leq l}p^{-1}(x_{i})\to 
X_{N}\setminus\bigcup_{1\leq i\leq l}\{x_{i}\}
$$ 
is biholomorphic, see [C], [R]. By definition, the domain 
$X_{M}:=p(M)\subset X_{N}$ is strongly pseudoconvex, and so it is
Stein. Without loss of generality we may assume that 
$x_{1},\dots, x_{l}\in X_{M}$. Thus 
$\cup_{1\leq i\leq l}\ p^{-1}(x_{i})=C_{M}$. 

Next, we introduce a complete K\"{a}hler
metric on the complex manifold $M\setminus C_{M}$ as follows.

First, according to [N] there is a proper one-to-one map 
$i:X_{M}\hookrightarrow\Co^{2n+1}$, $n=dim_{\Co}X_{M}$, which is an
embedding in regular points of $X_{M}$. Then $i(X_{M})$ is a complex
subvariety of $\Co^{2n+1}$. By $\omega_{e}$ we denote
the $(1,1)$-form on $M$ obtained as the pullback by 
$i\circ p$ of the Euclidean K\"{a}hler form on $\Co^{2n+1}$. Clearly,
$\omega_{e}$ is $d$-closed and positive outside $C_{M}$.

Similarly we can embed $X_{N}$ into $\Co^{2n+1}$ as a closed
complex subvariety. Let $j:X_{N}\hookrightarrow\Co^{2n+1}$ be an
embedding such that $j(X_{M})$ belongs to the open Euclidean ball $B$ of 
radius $1/4$ centered at $0\in\Co^{2n+1}$. 
Set $z_{i}:=j(x_{i})$, $1\leq i\leq l$. By
$\omega_{i}$ we denote the restriction to $M\setminus C_{M}$ of
the pullback with respect to $j\circ p$ of the form
$-\sqrt{-1}\cdot\partial\overline\partial\log(\log ||z-z_{i}||^{2})^{2}$
on $\Co^{2n+1}\setminus\{z_{i}\}$. (Here $||\cdot||$ stands for
the Euclidean norm on $\Co^{2n+1}$.) Since $j(X_{M})\subset B$, 
the form $\omega_{i}$ is K\"{a}hler.
Its positivity follows from the fact that the 
function $-\log(\log||z||^{2})^{2}$ is 
strictly plurisubharmonic for $||z||<1$.
Also, $\omega_{i}$ is extended to a smooth form on $M\setminus p^{-1}(x_{i})$.
Now, let us introduce a K\"{a}hler
form $\omega_{M}$ on $M\setminus C_{M}$ by the formula
\begin{equation}\label{e2}
\omega_{M}:=\omega_{e}+\sum_{1\leq i\leq l}\omega_{i}\ .
\end{equation}
\begin{Proposition}\label{p1}
The path metric $d$ on $M\setminus C_{M}$ induced by $\omega_{M}$ is complete.
\end{Proposition}
{\bf Proof.}
Assume, on the contrary, that there is a sequence $\{w_{j}\}$
convergent either to $C_{M}$ or to the boundary $bM$ of $M$ such that the
sequence $\{d(o,w_{j})\}$ is bounded (for a fixed point 
$o\in M\setminus C_{M}$). Then, since $\omega_{L}\geq\omega_{e}$, 
the sequence $\{i(p(w_{j}))\}\subset\Co^{2n+1}$ is
bounded. This implies that $\{w_{j}\}$ converges to $C_{M}$.
But since $\omega_{L}\geq\sum\omega_{i}$, the latter is impossible. One
can check this using single blow-ups of $\Co^{2n+1}$ at points $z_{i}$ and
rewriting the pullbacks to the resulting manifold
of $(1,1)$-forms  \penalty-10000
$-\sqrt{-1}\cdot\partial\overline\partial\log(\log ||z-z_{i}||^{2})^{2}$ 
in local coordinates
near exceptional divisors,
see, e.g., [GM] for similar arguments.\ \ \ \ \ $\Box$

Similarly one obtains complete K\"{a}hler metrics on unbranched 
coverings of $M\setminus C_{M}$ induced by
pullbacks to these coverings of the K\"{a}hler form $\omega_{M}$ 
on $M\setminus C_{M}$.\\
{\bf 3.3.} We retain the notation of the previous section. 

Let $r:N'\to N$ be an unbranched covering. Consider the corresponding
covering $(M\setminus C_{M})':=r^{-1}(M\setminus C_{M})$ of 
$M\setminus C_{M}$. 
We equip $(M\setminus C_{M})'$ with the complete
K\"{a}hler metric induced by the form
$\omega_{M}':=r^{*}\omega_{M}$. 
Next we consider the function $f:=\sum_{0\leq s\leq l}f_{s}$ 
on $(M\setminus C_{M})'$ such that
$f_{0}$ is the pullback by $i\circ p\circ r$ of 
the function $||z||^{2}$ on $\Co^{2n+1}$ and 
$f_{s}$ is the pullback by
$j\circ p\circ r$ of the function $-\log(\log ||z-z_{s}||^{2})^{2}$ on
$\Co^{2n+1}\setminus\{z_{s}\}$, $1\leq s\leq l$.
Clearly we have
\begin{equation}\label{e6}
\omega_{M}':=\sqrt{-1}\cdot\partial\overline\partial f\ .
\end{equation}

Let $E:=(M\setminus C_{M})'\times\Co$ be the trivial holomorphic line bundle
on $(M\setminus C_{M})'$. Let $g$ be the pullback to $(M\setminus C_{M})'$ of
a smooth plurisubharmonic function on $M$.
We equip $E$ with the hermitian metric $e^{f+g}$
(i.e., for $z\times v\in E$ the square of its norm in this
metric equals $e^{f(z)+g(z)}|v|^{2}$ where $|v|$ is the modulus of 
$v\in\Co$). Then the curvature $\Theta_{E}$ of $E$ satisfies
\begin{equation}\label{e7}
\Theta_{E}:=-\sqrt{-1}\cdot\partial\overline\partial
\log(e^{f+g})
=-\omega_{M}'-\sqrt{-1}\cdot
\partial\overline\partial g\leq -\omega_{M}'.
\end{equation}
From here by (\ref{e1'}) we obtain
\begin{equation}\label{e9}
H_{(2)}^{0,r}((M\setminus C_{M})',E)
=0\ \ \ {\rm for}\ \ \ r<n\ .
\end{equation}
{\bf 3.4.} In the proof we also use the following result.
\begin{Lm}\label{le1'}
Let $h$ be a nonnegative piecewise continuous function on $M$ equals 0 in 
some neighbourhood of $C_{M}$ and bounded on every compact subset of
$M\setminus C_{M}$. Then there exists a smooth
plurisubharmonic function $\hat g$ on $M$ such that 
$$
\hat g(z)\geq h(z)\ \ \ {\rm for\ all}\ \ \ z\in M.
$$
\end{Lm}
{\bf Proof.} Without loss of generality we identify $M\setminus C_{M}$
with $X_{M}\setminus\cup_{1\leq j\leq l}\ \{x_{j}\}$. Also, we identify
$X_{M}$ with a closed subvariety of $\Co^{2n+1}$ as in Section 3.2. 
Let $U$ be a neighbourhood of $\cup_{1\leq j\leq l}\ \{x_{j}\}$ such that
$h|_{U}\equiv 0$. By $\Delta_{r}\subset\Co^{2n+1}$ we denote the open
polydisk 
of radius $r$ centered at $0\in\Co^{2n+1}$. Assume without loss of generality
that $0\in X_{M}\setminus U$. Consider the monotonically 
increasing function
\begin{equation}\label{e15}
v(r):=\sup_{\Delta_{r}\cap X_{M}}h\ ,\ \ \
r\geq 0\ .
\end{equation}
By $v_{1}$ we denote a smooth monotonically increasing function satisfying
$v_{1}\geq v$ (such $v_{1}$ can be easily constructed
by $v$). Let us determine 
$$
v_{2}(r):=\int_{0}^{r+1}2v_{1}(2t)\ \!dt\ ,\ \ \ r\geq 0\ .
$$
By the definition $v_{2}$ is smooth, convex and monotonically
increasing. Moreover, 
$$
v_{2}(r)\geq\int_{\frac{r+1}{2}}^{r+1}2v_{1}(2t)\ \!dt\geq (r+1)v(r+1)\ .
$$
Next we define a smooth 
plurisubharmonic function $v_{3}$ on $\Co^{2n+1}$ by the formula
$$
v_{3}(z_{1},\dots, z_{2n+1}):=\sum_{j=1}^{2n+1}v_{2}(|z_{j}|)\ .
$$
Then the pullback of $v_{3}$ to $M$ is a smooth 
plurisubharmonic function on $M$. This is the required function 
$\hat g$. Indeed, under
the identification described at the beginning of the proof for
$|z|_{\infty}:=\max_{1\leq i\leq 2n+1}|z_{i}|$ we have
$$
\hat g(z)=v_{3}(z)\geq(|z|_{\infty}+1)v(|z|_{\infty}+1)\geq
\sup_{\Delta_{|z|_{\infty}+1}\cap X_{M}}h\geq h(z)\ \ \ {\rm for\ all}\ \ \ z\in M.
\ \ \ \ \ \Box
$$
{\bf 3.5.} In the proof 
of Theorem \ref{te2} we will assume without loss of generality
that $C_{M}$ is a {\em divisor with normal crossings}. 
Indeed, according to the Hironaka theorem, 
there is a {\em modification} $m: N_{H}\to N$ of $N$ from Section 1.3 
such that 
$m^{-1}(C_{M})$ is a divisor with normal crossings and
$m: N_{H}\setminus m^{-1}(C_{M})\to N\setminus C_{M}$ is biholomorphic.
By the definition $M_{H}:=m^{-1}(M)\subset N_{H}$ is strongly pseudoconvex.
Further, since $M$ is a complex manifold, $m$ induces an isomorphism of 
fundamental groups $m_{*}:\pi_{1}(M_{H})\to\pi_{1}(M)$. 
Thus for an unbranched covering $r:M'\to M$ of $M$ there are a covering $r_{H}:M_{H}'\to M_{H}$ and
a modification $m':M_{H}'\to M'$ such that $r\circ m'=m\circ r_{H}$ and
$m'$ induces an isomorphism of the corresponding fundamental groups.

Assume now that a $(0,1)$-form $\eta\in {\cal E}_{2}(M')$ satisfies the 
hypotheses of Theorem \ref{te2}. Consider its pullback 
$\widetilde\eta:=(m')^{*}\eta$ on $M_{H}'$. Clearly,
$\widetilde\eta$ also satisfies the hypotheses of Theorem \ref{te2} with $M$ 
replaced by $M_{H}$. Now, suppose that Theorem \ref{te2} is valid for $M_{H}'$, i.e., there is a continuous
function $\widetilde f\in {\cal L}_{2}(M_{H}')$ such that 
$\overline\partial\widetilde f=\widetilde\eta$ and $\widetilde f$
vanishes in a neighbourhood of $bM_{H}'$. Since by the definition of
$\eta$ the function $\widetilde f$ is holomorphic in a neighbourhood of
$(r\circ m')^{-1}(C_{M})\subset M_{H}'$ and $m':M_{H'}\to M'$ is a 
modification of $M'$, there is a function $f\in {\cal L}_{2}(M')$ such that
$\widetilde f=(m')^{*}f$. Obviously, 
$f$ satisfies the required statements of the theorem.\\
{\bf 3.6.} Let $U_{q}\subset\subset M$ be a simply connected
coordinate chart of $q\in C_{M}$
with complex coordinates $z=(z_{1},\dots,z_{n})$, $n=dim_{\Co}M$, such that
$z_{1}(q)=\cdots =z_{n}(q)=0$ and
\begin{equation}\label{normal}
C_{M}\cap U_{q}=\{f_{q}(z)=0\},\ \ \ f_{q}(z):=z_{1}\cdots z_{k}.
\end{equation}
(Such coordinates exist by the definition of a 
divisor with normal crossings.)

Let $\hat f$ be a function on $M\setminus C_{M}$ such that $r^{*}\hat f=f$, see Section 3.3. From the definition of $f$ we obtain
\begin{Lm}\label{le3.3}
$e^{\hat f}$ extended by 0 to $C_{M}$ is a continuous function on $M$ such that
$e^{\hat f}/|f_{q}|$ is unbounded on 
$U_{q}\setminus C_{M}$.\ \ \ \ \ 
$\Box$
\end{Lm}

Let $\omega$ be the associated (1,1)-form of a hermitian metric $g_{N}$ on $N$.
Since by the definition $\omega_{M}\geq\omega_{e}$ and the latter form
vanishes on $C_{M}$, we have locally near $C_{M}\cap U_{q}$
\begin{equation}\label{e3.8}
\omega_{M}^{n}\geq c'|f_{q}|^{2m'}\omega^{n}
\end{equation}
for some $c'>0$, $m'\in\N$. 
This and Lemma \ref{le3.3} imply that locally near $C_{M}\cap U_{q}$
\begin{equation}\label{e3.9}
e^{\hat f}\omega_{M}^{n}\geq c|f_{q}|^{2m}\omega^{n}
\end{equation}
for some $c>0$, $m\in\N$.

By $E_{n}(M)$ we denote a holomorphic line vector bundle on $M$ determined 
by the divisor $nC_{M}$,
$n\in\N$. Let $s_{1}$ be a holomorphic section of $E_{1}(M)$ 
defined in local coordinates on $U_{q}$ by functions $f_{q}$ from 
(\ref{normal}). Then $(r^{*}s_{1})^{n}$ is
a holomorphic section of the bundle $E_{n}'(M):=r^{*}E_{n}(M)$ on $M'$.

Since the hermitian bundle $E$ from Section 3.3 is
holomorphically trivial, we naturally identify sections of 
$E$ with functions on $(M\setminus C_{M})'$. Here and below we set
$X':=r^{-1}(X)$
for $X\subset M$. Also, the Banach space ${\cal L}_{2}(X')$ of continuous 
functions on $X'$ is defined similarly to ${\cal L}_{2}(M')$, see Section 1.3.

Let $(U_{i})_{i\in I}$ be a finite open cover of a neighbourhood
$\overline{M} (\subset\subset N)$ by simply connected coordinate charts
$U_{i}\subset\subset N$. 
\begin{Proposition}\label{pr3}
Suppose $h\in L_{2}((M\setminus C_{M})',E)$ is such that for any 
$U_{i}'$ there is a continuous function $h_{i}\in {\cal L}_{2}(U_{i}')$ 
so that $c_{i}:=h-h_{i}\in {\cal O}((U_{i}\setminus C_{M})')$.
Then there is an integer $n\in\N$ independent of $h$ such that 
$h\cdot (r^{*}s_{1})^{n}$ 
admits an extension
$\hat h\in C(M',E_{n}'(M))$. Moreover, $h|_{O'}\in {\cal L}_{2}(O')$ for
every $O\subset\subset M\setminus C_{M}$.
\end{Proposition}
{\bf Proof.} 
Let $U_{q}$ be a simply connected coordinate chart of $q\in C_{M}$
with the local coordinates satisfying (\ref{normal}).
We naturally identify $U_{q}'$ with
$U_{q}\times S$ where $S$ is the fibre of $r$. Then the hypotheses of the
proposition imply that
\begin{equation}\label{e3.10}
\int_{z\in U_{q}\setminus C_{M}}
\left(\sum_{s\in S}|h(z,s)|^{2}\right)e^{\hat f(z)+\hat g(z)}\ \!
\omega_{M}^{n}(z)<\infty
\end{equation}
where $\hat g$ is a smooth plurisubharmonic function on $M$ such that
$r^{*}\hat g=g$.\\
Diminishing if necessary $U_{q}$ assume that 
$\hat f$, $\omega_{M}^{n}$ satisfy (\ref{e3.9}) there.
Also, on $U_{q}$ we clearly have $\hat g\sim 1$. From here and 
(\ref{e3.10}) we obtain on $U_{q}$ 
\begin{equation}\label{e17}
\int_{z\in U_{q}\setminus C_{M}}
\left(\sum_{s\in S}|h(z,s)|^{2}\right)\ \!
|f_{q}(z)|^{2m}\omega^{n}(z)<\infty\ 
\end{equation}
where $f_{q}$ is defined by (\ref{normal}).

Further, according to the hypothesis of the proposition, there is a
continuous function $h_{q}\in {\cal L}_{2}(U_{q}')$ such that
$c_{q}:=h-h_{q}\in {\cal O}((U_{q}\setminus C_{M})')$. This and
(\ref{e17}) imply that every $f_{q}^{m}\cdot c_{q}(\cdot,s)$, $s\in S$,
is $L^{2}$ integrable with respect to the volume form
$(\sqrt{-1})^{n}\wedge_{i=1}^{n}dz_{i}\wedge d\overline{z}_{i}$. Using these
facts and the Cauchy integral formulas for coefficients of the Laurent
expansion of $f_{q}^{m}c_{q}(\cdot,s)$, one obtains easily that every
$f_{q}^{m}c_{q}(\cdot,s)$ can be extended holomorphically to $U_{q}$.
In turn, this gives a continuous extension $\hat h$ of 
$h\cdot (r^{*}f_{q})^{m}$ to $U_{q}'$.

Let $V_{q}\subset\subset U_{q}$ be another connected neighbourhood of
$q$. 
By the Bergman inequality for holomorphic functions, see, e.g., [GR, Chapter 6, Theorem 1.3], we have
\begin{equation}\label{ber2}
|h(y,s)f_{q}^{m}(y)|^{2}\leq A\int_{z\in U_{q}}|h(z,s)f_{q}^{m}(z)|^{2}
\omega^{n}(z)\ \ \ {\rm for\ all}\ \ \ (y,s)\in W_{q}'
\end{equation}
with some constant $A$ depending on $U_{q}$, $W_{q}$ and $\omega$ only.
Therefore from (\ref{e17}) and (\ref{ber2})
we obtain
$$
\sup_{z\in V_{q}}\left(\sum_{s\in S}|\hat h(z,s)|^{2}\right)^{1/2}<\infty.
$$
Equivalently, $\hat h|_{V_{q}'}\in {\cal L}_{2}(V_{q}')$.

Next assume that $U_{q}\subset (U_{i})_{i\in I}$ is a  
simply connected coordinate neighbourhood of a point 
$q\in M\setminus C_{M}$.
 Without loss of generality we may assume
that all such $U_{q}$ are relatively compact in $M\setminus C_{M}$.
Identifying
$U_{q}'$ with $U_{q}\times S$ we have anew inequality 
of type (\ref{e3.10}) for $h|_{U_{q}'}$. 
Since  $U_{q}\subset\subset M\setminus C_{M}$ and
$\hat f$, $\hat g$ and $\omega_{M}^{n}$ are smooth on $M\setminus C_{M}$ by their definitions, we obviously have on $U_{q}$
$$
e^{\hat f+\hat g}\cdot\omega_{M}^{n}\sim\omega^{n}\ .
$$
Similarly to (\ref{e17})-(\ref{ber2}) (with $f_{q}=1$) this implies that 
$h|_{V_{q}'}\in {\cal L}_{2}(V_{q}')$  for any connected neighbourhood 
$V_{q}\subset\subset U_{q}$ of $q$. 
Choose the above neighbourhoods $V_{q}$ so that they form a finite cover of 
a set $O\subset\subset M\setminus C_{M}$. Then from the implications
$h|_{V_{q}'}\in {\cal L}_{2}(V_{q}')$ we obtain
that $h|_{O'}\in {\cal L}_{2}(O')$. Now, choosing the
neighbourhoods $V_{q}$, $q\in C_{M}$, 
so that they form a finite cover of $C_{M}$ and taking as the 
$n$ the maximum of the numbers $m$ in the powers of $f_{q}$, see (\ref{e3.9}), 
we obtain that $h\cdot (r^{*}s_{1})^{n}$ admits an
extension $\hat h\in C(M',E_{n}'(M))$. By our construction
$n$ is independent of $h$.\ \ \ \ \ $\Box$\\
{\bf 3.7. Proof of Theorem \ref{te2}.}
Assume that a (0,1)-form $\eta$ belongs to 
${\cal E}_{2}(M')$, is $\overline\partial$-closed and
$r(supp\ \eta)\subset O\subset\subset M\setminus C_{M}$.

Let us define the function $g$ in the definition of the bundle
$E$ from Section 3.3 by Lemma \ref{le1'}. Namely, fix a neigbourhood 
$U\subset\subset M$ of $C_{M}$ and consider the function $h$ on $M$ defined by 
the formula
\begin{equation}\label{e3.81}
h(z):=\frac{\chi_{U^{c}}(z)}{dist(z,bM)}
\end{equation}
where $\chi_{U^{c}}$ is the characteristic function of $U^{c}:=M\setminus U$ and the distance to the boundary is defined by the path metric $d_{N}$
on $N$ induced by the Riemannian metric $g_{N}$. Further, according to Lemma \ref{le1'} we can find a $C^{\infty}$  plurisubharmonic function $\hat g$ on $M$ such that
$\hat g(z)\geq h(z)$ for all $z\in M$. Then in the definition of the metric on 
$E$ we choose $g:=r^{*}\hat g$.
\begin{Lm}\label{cons}
The form $\eta$ belongs to $L_{2}^{0,1}((M\setminus C_{M})',E)$.
\end{Lm}
{\bf Proof.}
We retain the notation of Proposition \ref{pr3}. Consider the set
$U_{q}'\cong U_{q}\times S$ on $M'$ for some $q\in M$ such that
$U_{q}\subset\subset M\setminus C_{M}$. Since $\eta\in {\cal E}_{2}(M')$, 
$r(supp\ \eta)\subset O\subset\subset M\setminus C_{M}$ and $\hat g$,
$\hat f$ and $\omega_{M}^{n}$ are bounded on $O$, for every such $U_{q}$
we have
\begin{equation}\label{cons1}
\int_{z\in U_{q}\setminus C_{M}}\left(\sum_{s\in S}|\eta|_{(z,s)}^{2}
\right)e^{\hat f(z)+\hat g(z)}\omega_{M}^{n}(z)<\infty.
\end{equation}
(Recall that $|\eta|_{(z,s)}^{2}$ stands for the norm of $\eta$ at
$(z,s)\in M'$ defined by the natural hermitian metric on the fibres of the
cotangent bundle $T^{*}M'$ on $M'$.) Taking a finite open cover of $O$ by
such sets $U_{q}$ we get the required statement.\ \ \ \ \ $\Box$

From Lemma \ref{cons} and the fact that $\overline\partial\eta=0$  we obtain 
by (\ref{e9}) that there exists a
function $F'\in L_{2}((M\setminus C_{M})',E)$ such that $\overline\partial F'=\eta$. Moreover, by the definition of $\eta$, this function is holomorphic on $(M\setminus C_{M})'\setminus r^{-1}(\overline{O})$. 
Also, since $\eta\in {\cal E}_{2}(M')$ the equation 
$\overline\partial G=\eta$ has local (continuous) solutions 
$f_{U}\in {\cal L}_{2}(U')$ for every $U\subset\subset M$ biholomorphic to
an open Euclidean ball of $\Co^{n}$. (In fact, since
$U'\cong U\times S$,
we can rewrite the equation $\overline\partial G=\eta$ on $U'$ as a 
$\overline\partial$-equation on $U$ with a measurable
Hilbert valued $(0,1)$-form on 
the right-hand side. Then we apply the formula presented in the proof of
Lemma 3.4 of [Br3] (see also [H, Section 4.2]) to solve this equation and to
get a solution from ${\cal L}_{2}(U')$, for similar arguments see 
[Br1, Appendix A].) 

Let us prove now
\begin{Lm}\label{zero}
There is a neighbourhood $U\subset M$ of $bM$ such that $F'|_{r^{-1}(U)}=0$.
\end{Lm}
{\bf Proof.} Let $q\in bM$ and $U_{q}\subset\subset N\setminus C_{M}$ be a simply connected coordinate chart of $q$.
Since $\pi_{1}(M)=\pi_{1}(N)$ by our assumption, the covering $M'$ of $M$ is contained in the corresponding covering $r:N'\to N$ of $N$. Thus $r^{-1}(U_{q})\subset N'$ can be naturally identified with $U_{q}\times S$ where $S$ is the fibre of $r$. Further, without loss of generality we may identify $U_{q}$ with an open Euclidean ball $B$ in $\Co^{n}$. In this identification, on each component $U_{q}\times\{s\}$, $s\in S$, the path metric $d$ on $N'$ is equivalent to the Euclidean metric on $B$ with the constants of equivalence independent of $s$. 

Next, for some $s\in S$ let us consider the restriction $F_{s}'$ of $F'$ to $U_{q}\times\{s\}=B$. We set $M_{s}':=M'\cap (U_{q}\times\{s\})$ and
$bM_{s}':=bM'\cap (U_{q}\times\{s\})$. 
Diminishing if necessary $U_{q}$, 
without loss of generality we may assume that these sets are connected.
Also by $dv$ we denote the Euclidean volume form on $\Co^{n}$. By the constructions of $\omega_{M'}$, see Section 3.2, and $f$, see Section 3.3, 
we clearly have
\begin{equation}\label{e316}
f|_{U_{q}\times\{s\}}\geq c\ \ \ {\rm and}\ \ \ \omega_{M'}^{n}|_{U_{q}\times\{s\}}\geq c\ \!dv
\end{equation}
for some $c>0$ independent of $s\in S$.

Further, by the definition $F_{s}'\in L_{2}(M_{s}',E)$. 
So by the choice of $g$ 
in the definition of the hermitian metric on $E$ using (\ref{e316}) we obtain
\begin{equation}\label{e317}
\int_{z\in M_{s}'}|F_{s}'(z)|^{2}e^{\frac{1}{dist(z,bM_{s}')}}dv(z)<\infty.
\end{equation}
Without loss of generality we may assume that $U_{q}\cap r(supp\ \eta)=\emptyset$. Thus $F_{s}'$ is holomorphic on $M_{s}'$ for each $s\in S$.
Now, from (\ref{e317}) using the mean-value property for the 
plurisubharmonic function $|F_{s}'|^{2}$ defined on $M_{s}'$ we easily
obtain that for any $y\in bM_{s}'$
\begin{equation}\label{zer}
\lim_{z\to y}F_{s}'(z)=0\ .
\end{equation}
Indeed, for a point $z$ sufficiently close to $y\in bM_{s}'$ consider a 
Euclidean ball
$B_{z}$ centered at $z$ of radius $r_{z}:=dist(z,bM_{s}')/2$. Choosing
$z$ closer to $y$ we may assume that $B_{z}\subset\subset M_{s}'$. Then
by the triangle inequality for the metric $d_{N}$ we have
$$
dist(w,bM_{s}')\leq 3r_{z}/2\ \ \ {\rm for\ all}\ \ \ w\in B_{z}.
$$
Now from (\ref{e317}) by the mean-value property
we get for some $c_{n}>0$ depending on $n$ only:
$$
\begin{array}{c}
\displaystyle
c_{n}r_{z}^{2n}e^{2/(3r_{z})}|F_{s}'(z)|^{2}\leq 
e^{2/(3r_{z})}\int_{w\in B_{z}}|F_{s}'(w)|^{2}dv(w)\leq\\
\\
\displaystyle
\int_{w\in B_{z}}|F_{s}'(w)|^{2}e^{\frac{1}{dist(w,bM_{s}')}}dv(w)\leq A<
\infty.
\end{array}
$$
Hence, 
$$
\lim_{z\to y} |F_{s}'(z)|^{2}\leq\lim_{z\to y}
\frac{Ae^{-2/(3r_{z})}}{c_{n}r_{z}^{2n}}=0.
$$
Thus (\ref{zer}) is true for any $y\in bM_{s}'$.

Next, since $M_{s}'$ is connected, (\ref{zer}) implies that 
$F_{s}'\equiv 0$ on 
$M_{s}'$ for each $s\in S$. Actually, let $z\in bM_{s}'$. Consider a complex
line $l_{z}$ passing through $z$ and containing the normal to $bM_{s}'$ at
$z$ (recall that $bM_{s}'$ is smooth). Then $l_{z}$ intersects $bM_{s}'$ 
transversely in a neighbourhood of $z$ in $bM_{s}'$.
This implies that there is a simply connected domain 
$W_{z}\subset l_{s}\cap M_{s}'$ whose boundary $bW$ contains $z$ such that
$F_{s}'|_{\overline{W}_{z}}\in C(\overline{W}_{z})$ and it
equals $0$ on an open subset of $bW_{z}$. Thus by the uniqueness
property for univariate holomorphic functions we have $F_{s}'=0$ on $W_{z}$.
Observe that if $z$ varies along $bM_{s}'$ the union of the connected 
components of $l_{z}\cap M_{s}'$ containing $W_{z}$ contains an open subset 
of $M_{s}'$. This implies that $F_{s}'\equiv 0$ on $M_{s}'$. 

Finally, taking a finite open cover of $bM$ by the
above sets $U_{q}$ and using similar arguments we obtain the required
neighbourhood $U$ of $bM$ (as the union of such $U_{q}$ intersected with $M$).
This completes the proof of the lemma.
\ \ \ \ \ $\Box$

Let us finish the proof of the theorem.
As established above, the function $F'$ 
satisfies conditions of Proposition \ref{pr3}. According to this
proposition there is a number $n\in\N$ independent of $F'$ such that
$F'\cdot (r^{*}s_{1})^{n}$
is extended to a continuous section of $E_{n}'(M)$ 
equals 0 on $r^{-1}(U)$. Moreover, $F'|_{O}\in {\cal L}_{2}(O')$
for any $O\subset\subset M\setminus C_{M}$.

We set 
$$
\widetilde F:=e^{F'}-1\ \ \ {\rm and}\ \ \ \widetilde\eta:=
\overline\partial\widetilde F=\widetilde F\eta.
$$
By the definitions of
$\eta$ and $F'$ we have $supp\ \widetilde\eta=supp\ \eta$ and
$\widetilde F$ is bounded on $supp\ \eta$. In particular,
$\widetilde\eta$ is $\overline\partial$-closed and belongs to 
$L_{2}^{0,1}((M\setminus C_{M})',E)$, as well. Then by (\ref{e9}) 
there is a function $\widetilde F'\in L_{2}((M\setminus C_{M})',E)$
such that $\overline\partial\widetilde F'=\widetilde\eta$. Applying
to $\widetilde F'$ the same arguments as to $F'$ we conclude that
$\widetilde F'|_{r^{-1}(U)}\equiv 0$ for some neigbourhood $U\subset M$
of $bM$. Since by the definition $\widetilde F-\widetilde F'$ is
holomorphic on $(M\setminus C_{M})'$ and equals zero on a neighbourhood
of $bM'$, from the connectedness of $M'$ we get $\widetilde F=\widetilde F'$.
Also, as in the case of $F'$, 
$\widetilde F'\cdot (r^{*}s_{1})^{n}$
is extended to a continuous section of $E_{n}(M')$.

Let $q$ be a regular point of $C_{M}':=r^{-1}(C_{M})$. 
The above properties
of $F'$ and $\widetilde F$ imply that for suitable
complex coordinates $z=(z_{1},\dots, z_{n})$ in a neighbourhood $U_{q}$ of $q$ 
we have $C_{M}'\cap U_{q}=\{z_{1}=0\}$ and
$$
e^{F'(z)}=z_{1}^{-n}A(z),\ \ \ 
F'(z)=z_{1}^{-n}B(z),\ \ \ z\in U_{q}\setminus C_{M},
$$
where $A, B\in {\cal O}(U_{q})$.
Suppose that $A(z)=z_{1}^{l}A'(z)$ for some $0\leq l<n$ with 
$A'\in {\cal O}(U_{q})$ not identically $0$ on $C_{M}'\cap U_{q}$.
Then there is a point $p\in C_{M}'\cap U_{q}$ and its neighbourhood
$W\subset U_{q}$ so that $A'(z)\neq 0$ for all $z\in\overline{W}$.
Thus we can 
introduce complex coordinates $y=(y_{1},\dots, y_{n})$ on $W$ by 
$y_{1}:=z_{1}(A'(z))^{1/(n-l)}$,
$y_{2}=z_{2},\dots, y_{n}=z_{n}$. In these coordinates we have
$e^{F'(y)}=y_{1}^{l-n}$, $y\in W\setminus C_{M}$. 
Since $F'\in {\cal O}(W\setminus C_{M})$, the latter
is impossible. This contradiction shows that
$l\geq n$ and so $e^{F'}$ admits a holomorphic extension to $U_{q}$. From here
we obtain easily that $F'$ admits a holomorphic extension to $U_{q}$, as well.

Taking an open cover of regular points of $C_{M}'$ by such neighbourhoods
$U_{q}$, from the above arguments we obtain that $F'$ is extended  
holomorphically to $C_{M}'$ (it is extended to nonregular points of $C_{M}'$ 
by the Hartogs theorem because
the complex codimension of the set of such points in $M'$ is $\geq 2$).

Finally, the extended function $F$ (i.e., the extension of $F'$) belongs
to ${\cal L}_{2}(M')$. Indeed, by the definition 
$F|_{O'}\in {\cal L}_{2}(O')$ for every $O\subset\subset N\setminus C_{M}$.
Assume now that $q\in C_{M}$. Let $U$ be a simply connected
coordinate chart of $q$ with complex coordinates $z=(z_{1},\dots, z_{n})$
such that $z_{1}(q)=\cdots=z_{n}(q)=0$, $C_{M}\cap M=\{z_{1}\cdots z_{k}=0\}$
and $\overline{U}=\{z\in M\ :\ \max_{1\leq k\leq n}|z_{k}|\leq 1\}$. We
identify $\overline{U}$ with the unit polydisk in $\Co^{n}$ and by
$\To^{n}$ we denote its boundary torus. Also, we naturally identify
$(\overline{U})'\subset M'$ with $\overline{U}\times S$ where
$S$ is the fibre of $r:M'\to M$. Diminishing, if necessary, $U$ we will
assume that $F$ is holomorphic in a neighbourhood $O':=r^{-1}(O)$ of 
$(\overline{U})'$ where $O$ is a neighbourhood of $\overline{U}$.

Let $\{S_{l}\}_{l\in\N}\subset S$ be an
increasing sequence of finite subsets of $S$ such that $\cup_{l}S_{l}=S$.
Then from the Cauchy integral formula we obtain
$$
\lim_{l\to\infty}\left(\sum_{s\in S_{l}}|F(y,s)|^{2}
\right)\leq \left(\frac{1}{2\pi}\right)^{n}
\int_{x\in\To^{n}}\sum_{s\in S}\frac{|F(x,s)|^{2}}{(1-|z_{1}(y)|)\cdots 
(1-|z_{n}(y)|)}\ \!dx,\ \ \ y\in U,
$$
where $dx$ is the volume form on $\To^{n}$. Since $\To^{n}\subset\subset
M\setminus C_{M}$, $F|_{\To^{n}}\in {\cal L}_{2}(\To^{n}\times S)$.
This implies that $F|_{r^{-1}(y)}\in l_{2}(S)$ for all 
$y\in V_{q}:=\{z\in U\ :\ \max_{1\leq k\leq n}|z_{k}|\leq 1/2\}$ and 
the $l_{2}$ norms $|\cdot|_{y}$ of these functions are uniformly bounded.
Choosing a finite cover of $C_{M}$ by such $V_{q}$ and taking into account
that $F|_{O'}\in {\cal L}_{2}(O')$ for every 
$O\subset\subset N\setminus C_{M}$, from the above we obtain
that $F\in {\cal L}_{2}(M')$. Also, $\overline\partial F=\eta$.
 
This completes the proof of the theorem.
\ \ \ \ \ $\Box$

\end{document}